\documentclass[1p]{elsarticle}

\usepackage{amsmath,amssymb,amsfonts,amsthm,amsgen,amscd}
\usepackage[matrix,arrow,curve,frame]{xy}
\usepackage{graphicx}

\usepackage[USenglish]{babel}

\usepackage{mathtools}

\usepackage{cite}
\usepackage{enumitem} 

\usepackage{float}
\usepackage{microtype}

\usepackage{bbm}

\numberwithin{equation}{section}
\numberwithin{figure}{section}

	\newcommand{\Z}{\mathbb{Z}}
	\newcommand{\CP}{\mathbbm{C}\mathrm{P}}
	\newcommand{\C}{\mathbbm{C}}
	\newcommand{\RP}{\mathbbm{R}\mathrm{P}}
	\newcommand{\R}{\mathbbm{R}}
	\newcommand{\KP}{\mathbbm{K}\mathrm{P}}
	\newcommand{\K}{\mathbbm{K}}
	\renewcommand{\H}{\mathbbm{H}}
	\newcommand{\dis}{\displaystyle}
	\newcommand{\scr}{\scriptstyle}

	\newcommand{\incl}{\operatorname{incl}}
	\newcommand{\id}{\operatorname{id}}

	\newcommand{\MC}{\operatorname{MC}}
	\newcommand{\MCC}{\operatorname{MCC}}
	\newcommand{\Ker}{\operatorname{Ker}}
	\newcommand{\comp}{\, \scriptstyle \circ \displaystyle}

	\newtheorem{thm}{Theorem}[section]
	\newtheorem{cor}[thm]{Corollary}
	
	\newtheorem{prop}[thm]{Proposition}
	
  \theoremstyle{definition}
	\newtheorem{defn}[thm]{Definition}

	\newtheorem{rem}[thm]{Remark}
	\newtheorem{rem_cor}[thm]{Remark and Correction}
  \theoremstyle{remark}
	
  \newtheoremstyle{example}{3pt}{3pt}{}{}{\bfseries}{:}{.5em}{}
	\theoremstyle{example}
	\newtheorem{exa}[thm]{Example}
	
	\newtheorem*{example}{Example}
	\newtheorem*{question}{Question}
	
	\newcounter{countertitem}[thm]
	\newcommand{\titem}{%
		\stepcounter{countertitem}
		{\rm \Roman{countertitem}.}
	}

\usepackage{todonotes}

\begin{document}

\begin{frontmatter}

\author[enc]{Ulrich Koschorke}
\title{Nielsen numbers in topological coincidence theory}

\address[enc]{Department Mathematik, Universit\"at Siegen,
    57068 Siegen, Germany}

\ead{koschorke@mathematik.uni-siegen.de}

\begin{abstract}
  We discuss coincidences of pairs
  $ (f_1, f_2) $
  of maps between manifolds.
  We recall briefly the definition of four types of Nielsen numbers which arise naturally from the geometry of generic coincidences.
  They are lower bounds for the minimum numbers MCC and MC which measure to some extend the 'essential' size of a coincidence phenomenon.
  
  In the setting of fixed point theory these Nielsen numbers all coincide with the classical notion but in general they are distinct invariants.
  
  We illustrate this by many examples involving maps from spheres to the real, complex or quaternionic projective space
  $ \KP(n') $.
  In particular, when
  $ n' $
  is odd and
  $\, \K \,=\, \R \,$
  or
  $\, \C \,$,
  or when
  $ n' \equiv 23 \mod 24 $
  and
  $ \K \,=\; \H \,$,
  we compute the minimum number MCC and all four Nielsen numbers for every pair of these maps, and we establish a 'Wecken theorem' in this context
  (in the process we correct also a mistake in previous work concerning the quaternionic case).
  However, when
  $ n' $
  is even, counterexamples can occur, detected e.g. by Kervaire invariants.
\end{abstract}


\end{frontmatter}

\section{Introduction and discussion of results}

Throughout this paper let
$ f_1, f_2 : M^m \longrightarrow N^n $
be (continous) maps between connected smooth manifolds (of the indicated dimensions
$ m, n \geq 1 $)
without boundary, 
$ M $
being compact.

Consider the {\itshape coincidence set}
\stepcounter{thm}
\begin{equation}\label{eqn:1.1}
  C(f_1, f_2) {:=} \left\lbrace x \in M \vert f_1(x) =f_2(x) \right\rbrace.
\end{equation}
Its size and shape may vary greatly when we deform
$ f_1 $ and
$ f_2 $.
However, in topological coincidence theory we are not interested in any such 'inessential' changes.
We would like to capture those features which remain unchanged by arbitrary homotopies.
One possible measure of the size is the {\itshape{\bfseries m}inimum number of {\bfseries c}oincidence points}
\stepcounter{thm}
\begin{equation}\label{eqn:1.2}
  MC(f_1. f_2) = \min\{\# C(f_1',f_2') \,\vert\, f_1' \sim f_1, f_2' \sim f_2 \}.
\end{equation}
It follows from a result of R. Brooks [Br] 
that we obtain the same minimum number if we deform only one of the two maps
$ f_1, f_2 $
by a homotopy while leaving the other map fixed.

\begin{example} \textbf{fixed points.}	
  Let
  $ f $
  be a selfmap of
  $ M $.
  Then
  \begin{equation*}
    MC(f,\id) = MF(f) \,{:=}\, \min\big\{\#\{ x \in M\ \vert\, f'(x) = x \} \;\big|\; f' \sim f\, \big\}
  \end{equation*}
  is the classical minimum number of fixed points which plays a central role in topological fixed point theory (cf. e.g. [N], [Ji 1-3], [Ke], [Z] and [B1], p.9).  
  \hfill$\Box$
\end{example}

In coincidence theory we do not assume that the dimensions of
$ M $ and
$ N $
are equal.
Thus
$ MC(f_1, f_2) $
may often be infinite and hence a rather crude invariant (generically the coincidence set is an
$ (m-n) $--dimensional manifold!).
A sharper measure for essential coincidence phenomena seems to be the {\bfseries m}inimum number of {\bfseries c}oincidence (path--){\bfseries c}omponents
\stepcounter{thm}
\begin{equation}\label{eqn:1.3}
  MCC(f_1.f_2) = \min\{\#\pi_0(C(f_1',f_2')) \vert f_1' \sim f_1, f_2' \sim f_2 \}
\end{equation}
which is always finite (due to the compactness of the domain
$ M $).

These minimum numbers are the {\itshape principal object of study} in topological coincidence theory (compare [B1], p.9). 
The case when they vanish is of particular interest:
\begin{defn}\label{def:1.4}
  The pair
  $ (f_1. f_2) $
  of maps is called {\itshape loose} if there are homotopies
  $ f_1 \sim f_1', f_2 \sim f_2' $
  such that
  $ f_1'(x) \neq f_2'(x) \text{ for all } x \in M $
  (i.e.
  $ f_1, f_2 $
  can be 'deformed away' from one another).
\end{defn}

Just as in fixed point theory, the determination of minimum numbers can be helped greatly by a very natural decomposition of the coincidence set into 'Nielsen classes' and by a resulting notion of Nielsen numbers.
These are based on a careful geometric analysis of generic coincidence data, as follows (for more details see e.g. [K2], [K3]).	

After small approximations we may assume that both
$ f_1 $ and
$ f_2 $
are smooth and that the map
\begin{equation*}
  (f_1, f_2) : M \longrightarrow N \times N 
\end{equation*}
is transverse to the diagonal
\begin{equation*}
  \Delta = \{( y_1, y_2) \in N \times N \,\vert\, y_1 = y_2 \}.
\end{equation*}
Then
$ C(f_1, f_2) = (f_1, f_2)^{-1}(\Delta) $
is a smooth submanifold of M.
\addtocounter{figure}{4}
\begin{figure}[H]
  \centering
  \includegraphics[width=0.8\textwidth]{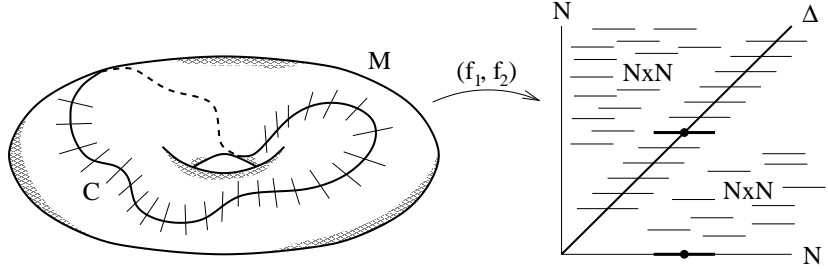}
  \caption{A generic coincidence manifold and its normal bundle}
  \label{fig:figure1}
\end{figure}

Our first coincidence datum keeps track of the smooth embedding
\addtocounter{equation}{2}
\addtocounter{thm}{2}
\begin{subequations}\label{equ:1.6}
  \renewcommand{\theequation}{\theparentequation,\roman{equation}}
  \begin{equation}\label{equ:1.6i}
    g : C(f_1, f_2) \hookrightarrow M.
  \end{equation}
  The normal bundle
  $ \nu(g) \text{ of } g $
  is described by the composite vector bundle isomorphism
  \begin{equation}\label{equ:1.6ii}
    \bar{g}^{\#} \;:\; \nu(g) \longrightarrow (f_1, f_2)^* (\nu(\Delta, N \times N)) \cong f_1^*(T N)
  \end{equation}
  induced by the tangent map of
  $ (f_1, f_2) $.
  
  Finally, there is a lifting
  \begin{equation}\label{equ:1.6iii}
    \widetilde{g} : C(f_1, f_2) \longrightarrow E(f_1, f_2)
  \end{equation}
  of $ g $, defined by
  $ \,\widetilde{g}(x) \;{:=}\; (x, $
  constant path at
  $ f_1(x) = f_2(x)\, $)\,; \,
  here	
  \begin{equation*}
    E(f_1, f_2) {:=} \{(x, \theta) \in M \times P(N) \,\vert\, \theta(0) = f_1(x), \theta(1) = f_2(x) \}
  \end{equation*}
  and
  $ P(N) $
  denotes the space of all continuous paths
  $ \theta : [0,1] \longrightarrow N $,
  with the compact--open topology.
  Though it may look innocuous, this third datum
  $\widetilde{g} $
  is by no means negligeable.
  It yields not only the Nielsen decomposition, but also important extra information (being responsible for the sometimes striking difference between the Nielsen numbers
  $ \widetilde{N}(f_1, f_2) $ and 
  $ N(f_1, f_2) $,
  cf. e.g. theorem \ref{thm:1.18}, corollary \ref{cor:1.24} and example \ref{exa:1.27} below).
\end{subequations}

The three data (\ref{equ:1.6}, i-iii) represent the nonstabilized normal bordism class
\stepcounter{thm}
\begin{equation}\label{equ:1.7}
  \omega^{\#}(f_1, f_2) = [C(f_1, f_2), \widetilde{g}, \bar{g}^{\#} ]\; \in\; \Omega^{\#}(f_1, f_2)	
\end{equation}
in a suitable bordism set
$ \Omega^{\#}(f_1, f_2) $
(for more details concerning this and the following constructions see [K3] and [K2]).	

If we keep track of
$ g $ and
$ \bar{g}^{\#} $
only as a continuous map and a {\itshape stable} vector bundle isomorphism we get the invariant
\stepcounter{thm}
\begin{equation}\label{equ:1.8}
  \widetilde{\omega}(f_1, f_2)\; \in\; \Omega_{m-n}(E(f_1, f_2);\widetilde{\varphi})	
\end{equation}
in a (standard) normal bordism group (with coefficients in a suitable virtual vector bundle
$ \widetilde{\varphi} $).

If we forget also the lifting
$ \widetilde{g} $
we obtain the normal bordism class
\stepcounter{thm}
\begin{equation}\label{equ:1.9}
  \omega(f_1, f_2) \;\in\; \Omega_{m-n}(M;\, \varphi = f_1^*(TN)-TM).
\end{equation}

Finally, by applying the Hurewicz homomorphism
$ \mu $
we may extract the invariant
\stepcounter{thm}
\begin{equation}\label{equ:1.10}
  \omega^{\Z}(f_1, f_2) = \mu(\omega(f_1, f_2)) \in H_{m-n}(M; \widetilde{\Z}_\varphi)
\end{equation}
in homology with integer coefficients (which are twisted like
$ \varphi $,
cf. \ref{equ:1.9}).

Each of these
$ \omega $--invariants
depends only on the homotopy classes of
$ f_1 $ and
$ f_2 $
and vanishes if the pair
$ (f_1, f_2) $
is loose (cf. definition \ref{def:1.4}).

Frequently the 'root' case where one of the maps
$ f_1, f_2 $
has a constant value
$* \in N $,
plays an important role.
\begin{defn}\label{def:1.11}
  Given a map
  $ f: M \longrightarrow N $,
  we define
  $ \deg^{\#}(f) = \omega^{\#}(f,*) $
  and similarly for
  $ \widetilde{\deg}(f),\, \deg(f) \text{ and } \deg^{\Z}(f) $.
\end{defn}
In the general case of arbitrary
$ f_1, f_2 $
the looseness obstruction
$ \omega^{\#}(f_1, f_2) $
contains often much more information than the other, increasingly weaker,
$ \omega $--invariants;
but it is also hardest to handle (in general
$ \Omega^{\#}(f_1, f_2) $
need not even be a group).
For the sake of simplification, let us extract numerical invariants (which will turn out to be useful bounds for minimum numbers).
\begin{defn}\label{def:1.12}
  The set
  $ \pi_0(E(f_1, f_2)) $
  of path components of the space
  $ E(f_1, f_2) $ (cf. \ref{equ:1.6iii})
  is called {\itshape Reidemeister set of the pair
  $ (f_1, f_2) $}.
  Its cardinality
  (in $\{0,1,\ldots,\infty\}$)
  is the {\itshape Reidemeister number}
  $ R(f_1, f_2) $.
\end{defn}

If 
$ x_0 \in M $
is a coincidence point put
$\, y_0 \,{:=}\, f_1(x_0) = f_2(x_0) \in N $.
According to [K2], 2.1, 
there exists a canonical bijection
\begin{equation*}
  \pi_1(N,y_0)/\text{Reidemeister equivalence }\; \longleftrightarrow\; \pi_0(E(f_1, f_2))	
\end{equation*}
where we call 
$ [\theta],[\theta'] $
{\itshape Reidemeister equivalent} if 
$ [\theta'] = f_{1*}(\gamma)^{-1} \cdot [\theta] \cdot f_{2*}(\gamma) $
for some
$ \gamma \in \pi_1(M,x_0)) $.
Thus \ref{def:1.12} gives just a base point free version of the standard definition of Reidemeister sets and numbers.

If 
$ M $
happens to be simply connected then the Reidemeister number depends only on the target manifold 
$ N $
and we have
\renewcommand{\theequation}{\arabic{section}.\arabic{equation}'}
\stepcounter{equation}
\begin{equation}\label{equ:1.12}
  R(f_1, f_2) \,\equiv\, R_N\; {:=}\; \# \pi_1(N).	
\end{equation}

Next we observe that the decomposition of
$ E(f_1, f_2) $
into its path components yields a disjoint decomposition of the coincidence set
$ C(f_1, f_2) $
into its parts
$ \widetilde{g}^{-1}(A), A\in \pi_0(E(f_1, f_2)) $.
In the generic case, these parts are closed
$ (m-n)$--submanifolds of
$ M $;
their (restricted) coincidence data as in (\ref{equ:1.6}, i-iii)
contribute to the
$ \omega$--invariants
defined in \eqref{equ:1.7}--\eqref{equ:1.10}.

\begin{defn}\label{def:1.13}
  {\itshape The Nielsen number}
  $\, N^{\#}(f_1, f_2) $
  (or
  $\, \widetilde{N}(f_1, f_2),\, N(f_1, f_2),\, N^{\Z}(f_1, f_2) \text{ }$,
  resp.) is the number of pathcomponents
  $ A \text{ of } E(f_1, f_2) $
  such that the contribution of
  $ \widetilde{g}^{-1}(A) \text{ to } \omega^{\#}(f_1, f_2) \text{ (or } \widetilde{\omega}(f_1, f_2), \omega(f_1, f_2), \omega^{\Z}(f_1, f_2) \text{, resp.)}$
  is nontrivial ('essential').
  
  {\bfseries Warning (change of notation).}	
  Until 2010 I denoted the Nielsen number
  $\widetilde{N}(f_1, f_2) $
  (which is based on 
  $ \widetilde{\omega}(f_1, f_2)) \text{ by } N(f_1, f_2) $.
\end{defn}

When 
$ m = n $
each of our four types of Nielsen numbers coincides with the classical notion of a Nielsen number which is so central e.g. in topological fixed point theory.

However, in strictly positive codimensions
$ m-n > 0 $,
we get four distinct types of Nielsen numbers which are lower bounds of the minimum and Reidemeister numbers (cf. [K3], 
theorem 1.2, and [K6]). 
Indeed,\vspace{-1cm}
\numberwithin{equation}{section}
\stepcounter{equation}
\stepcounter{thm}
\begin{equation}\label{equ:1.14}\vspace{-1cm}
  \begin{aligned}
    \xymatrix {
      & & \\
      & & \\
      MC \ar@<-5.5pt>@{}[r]^-{\overset{\displaystyle\not\equiv}{\displaystyle{\geq}}} &
      MCC \ar@<-5.5pt>@{}[r]^-{\overset{\displaystyle\not\equiv}{\displaystyle{\geq}}} \ar@{}[d]|{\textstyle{\text{if }n\neq 2:\,\text{\rotatebox{90}{$\geq$}}\qquad\qquad\;}} \ar@{}[dd]|{\textstyle{R}} \ar@{}[u]|{\rotatebox{90}{$<$}} \ar@{}[uu]|{\textstyle{\infty}} &
      N^{\#} \ar@<-5.5pt>@{}[r]^-{\overset{\displaystyle\not\equiv}{\displaystyle{\geq}}} &
      \widetilde{N} \ar@<-5.5pt>@{}[r]^-{\overset{\displaystyle\not\equiv}{\displaystyle{\geq}}} &
      N \ar@<-5.5pt>@{}[r]^-{\overset{\displaystyle\not\equiv}{\displaystyle{\geq}}} &
      N^{\Z} \ar@<-5.5pt>@{}[r]^-{\overset{\displaystyle\not\equiv}{\displaystyle{\geq}}} &
      0 \\
      & & \\
      & &
    }
  \end{aligned}
\end{equation}
where
$ N^{\Z} $
seems to vanish most of the time (except maybe when e.g. aspherical manifolds such as tori are involved).

This suggests a very natural {\bfseries two--step program} for investigating minimum numbers.
First we have to decide when
$ MCC(f_1, f_2) \text{ (or even } MC(f_1, f_2)\text{)} $
is equal to one of the Nielsen numbers and to which one (such results are costumarily called {\itshape 'Wecken theorems'} in honor of F. Wecken and his work, cf. [We]).	
Secondly, we must determine the relevant Nielsen number.
(Here it is helpful that the possible values of Nielsen numbers are often severely restricted).

\begin{exa}\label{exa:1.15} $\mathbf{ M = S^m, N = S^n, m,n \geq 1.}$
  Let a denote the antipodal involution on the sphere
  $ S^n$.
  Then
  \begin{equation*}
    MCC(f_1, f_2) = N^{\#}(f_1, f_2) = \begin{cases}
                                         1 & \text{ if } n \neq 1 \text{ and } f_1 \not\sim a \scr \circ \dis f_2; \\
                                         \left| d^{\scr \circ \dis}(f_1) - d^{\scr \circ \dis}(f_2) \right| & \text{ if } m = n = 1; \\
                                         0 & \text{ otherwise};
                                       \end{cases}
  \end{equation*}
  (here 
  $ d^{\scr \circ \dis}(f_i) \in \Z $
  denotes the usual degree).
  Moreover
  \begin{equation*}
    MC (f_1,f_2) =
      \begin{cases}
        0 & \text{if } f_1 \sim a \comp f_2; \\
        1 & \text{if } m,n \geq 2 \,\text{ and } [f] \in E (\pi_{m-1} (S^{n-1}))\backslash\{0\}; \\
        \vert d^0 (f_1) - d^0 (f_2) \vert & \text{if } m = n = 1; \\
        \infty & \text{if }m > n \geq 2 \text{ and } [f] \not \in E (\pi_{m-1}(S^{n-1}));
      \end{cases}
  \end{equation*}
  (here
  $ [f] \,{:=}\, [f'_1] - [a \comp f'_2] \, \in \, \pi_m(S^n) $
  where the basepoint preserving maps
  $ f'_1 $
  and
  $ a \scr \circ \dis f'_2 $
  are (freely) homotopic to
  $ f_1 $
  and
  $ a \comp f_2 $,
  resp.).
  
  If
  $ m = n $,
  then
  \begin{equation*}
    \MC(f_1, f_2) = \MCC(f_1, f_2) = N^{\#}(f_1, f_2) = \widetilde{N}(f_1, f_2) = N(f_1, f_2) = N^{\Z}(f_1, f_2).
  \end{equation*}
  
  On the other hand {\itshape assume} that
  $ (m,n) \neq (1, 1) $.
  Then we have:
  \begin{equation*}
    \widetilde{N}(f_1, f_2) = \begin{cases}
                                0 &\text{ if } \;\Gamma(f_1) = \Gamma(a \comp f_2); \\
                                1 &\text{ otherwise};
                              \end{cases}
  \end{equation*}
  (here
  \begin{equation*}
    \Gamma \,\, {:=}\, \oplus E^{\infty} \scr \circ \dis \gamma_{k} \;:\; [S^m, S^n] \cong \pi_m(S^n) \longrightarrow \underset{k \geq 1}{\oplus} \pi_{m-1-k(n-1)}^S
  \end{equation*}
  where $ E^{\infty} \comp \gamma_{k} $ denotes the stabilized $k^{\text{th}}$ Hopf--James invariant homomorphism);
  \begin{align*}
    N(f_1, f_2) &= \begin{cases}
                    0 &\text{ if } \; E^{\infty}([f_1]) = (-1)^{n+1} E^{\infty}([f_2]); \\
                    1 &\text{ otherwise};
                  \end{cases}\\
    N^{\Z}(f_1, f_2) &= 0 \,\text{ unless } m=n.
  \end{align*}
  
  This follows from [K2], 1.14	
  (see also [K6], 1.10).	
  \hfill$\Box$
\end{exa}

For further illustrations let us consider the more general case where
$ M = S^m$,
but no restrictions are put on
$ N $.
When 
$ m \text{ or } n \text{ equal } 1 $
then both minimum numbers MC and MCC as well as the four Nielsen numbers vanish identically, except in the case
$ M = N = S^1 $
where all these numbers are equal to
$ \left| d^{\scr \circ \dis}(f_1) - d^{\scr \circ \dis}(f_2) \right| $
for any selfmaps
$ f_1, f_2 $
of the circle
$ S^1 $
(compare example \ref{exa:1.15} above).

Thus we may assume that
$ m, n \geq 2 $
in further discussions.
Then 
$ S^m $
is simply connected and the Reidemeister number agrees with the order of
$ \pi_1(N) $
(cf. \ref{equ:1.12}).
Furthermore, given a triple
$ (C, \bar{g}^{\#}, \widetilde{g}) $
as in \ref{equ:1.6},i--iii, the
$ n$--codimensional
submanifold
$ C \text{ of } S^m $
allows a retraction
$ r $
(unique up to homotopy) to a point
$ x_0 \in S^m $.
Thus the choice of an orientation for the tangent space
$ T_{f_1(x_0)}(N) $
determines a trivialization
$ \bar{g}^{\#'} $
of the normal bundle
$ \nu(g) \text{ of } C \text{ in } S^m $
(cf. \ref{equ:1.6ii}).
Moreover the adjoint of the map
$ \widetilde{g} $,
suitably concatenated with the homotopies
$ f_1 \scr \circ \dis r \text{ and } f_2 \scr \circ \dis r $,
yields a map
$ \widetilde{g}' \text{ from } C $
into the loop space
$ \Omega N \text{ of } N $.
Then the bordism classes of the triples
$ (C, \bar{g}^{\#}, \widetilde{g}) \text{ and } (C, \bar{g}^{\#'}, \widetilde{g}')) $
determine one another.

As usual the {\bfseries Pontrjagin--Thom procedure} allows us to translate this geometric description of coincidence data into the language of homotopy theory.
Let 
$ C \times \R^n \subset S^m $
be a tubular neighborhood of
$ C = C \times \{0\} $
compatible with
$ \bar{g} $.
Also define a map
$ h \text{ from } S^m $
into the Thom space
$ ((\Omega N) \times \R^n) \cup \{\infty\} $
(of the trivial n--plane bundle over
$ \Omega N $)
by
\begin{equation*}
  h(x) {:=} \begin{cases}
              (\widetilde{g}'(c), v) & \text{if } x = (c, v) \in C \times \R^n; \\
              \infty & \text{if } x \;\in\; S^m \ \backslash\ C \times \R^n.
            \end{cases}
\end{equation*}
This Thom space can be identified with the smash product
$ S^n \wedge \Omega N^+ $
of the (pointed) spaces
$ S^n \text{ and } \Omega N^+ (=\Omega N, \text{ with an extra point } + \text{ added}) $.
Then the homotopy class
$ [h] \in \pi_m(S^n \wedge \Omega N^+) $
determines and is determined by the bordism class
$ [C, \bar{g}^{\#'}, \widetilde{g}'] $
or, equivalently, 
$ [C, \bar{g}^{\#}, \widetilde{g}] $.
For more details (also concerning base points) see e.g. proposition 2.5 (and the appendix) in [K3].	

Similarly, the (stabilized) invariants
$ \widetilde{\omega}(f_1, f_2) \text{ and } \omega(f_1, f_2) $---when
translated from the language of framed bordism groups to homotopy theory via the Pontrjagin--Thom procedure---take values in (stable) homotopy groups.
Then our four
$ \omega$--invariants
fit into the commuting diagram \ref{equ:1.16} of group homomorphisms (where
$ m, n \geq 2 $).
\stepcounter{equation}
\stepcounter{thm}
\begin{equation}\label{equ:1.16}
    \begin{aligned}
    \xy
      (47,70) *+++{\pi_m(S^n \wedge \Omega N^+)}="12" ;
      (60,50) *+++{{\lim\limits_{k \rightarrow \infty} \pi_{m+k}(S^{n+k} \wedge \Omega N^+) = \pi_{m-n}^S(\Omega N)}}="22" ;
      (0,40)  *+++{{\pi_m(N \times N)}}="31" ;
      (55,30) *+++{{\lim\limits_{k \rightarrow \infty} \pi_{m+k}(S^{n+k}) = \pi_{m-n}^S(\text{point})}}="42" ;
      (47,10) *+++{H_{m-n}(S^m; \Z) = }="52" ;
      (73,13) *+{\Z \;\;\text{  if } m = n;} ;
      (73,7) *+{0   \;\;\text{  if } m \neq n.} **\frm{\{};
      \ar @{->} "31";"12"^-{\omega^{\#}}
      \ar @{->} "31";"22"^-{\widetilde{\omega}}
      \ar @{->} "31";"42"_-{\omega}
      \ar @{->} "31";"52"_-{\omega^{\Z}}
      \POS (50,67) \ar @{->} (50,54)^-{E^{\infty}}
      \POS (50,47) \ar @{->} (50,34)^-{(\text{constant map})_*}
      \POS (50,27) \ar @{->} (50,14)^-{\mu}
    \endxy
  \end{aligned}
\end{equation}

When
$ m = n $
and
$ \pi_1(N) = 0 $ then
the vertical arrows in diagram \ref{equ:1.16} are isomorphisms and the four
$ \omega$--invariants
have equal strength.
However, when
$ m \neq n $
the (classical) homological looseness obstruction
$ \omega^{\Z} $
is completely useless; in contrast, the other
$ \omega$--invariants---and
in particular
$ \omega^{\#} $---allow
us often to compute the minimum number MCC.

\begin{exa}\label{exa:1.17}
  {\bfseries real, complex or quaternionic projective spaces.} Let
  $ M = S^m,\, N^n =\KP(n'), \,\K = \R, \,\C \text{ or } \H,\; m, n \geq 2. $
  Here
  $ n = d n' $
  where
  \begin{equation*}
    d {:=} \dim_{\R}(\K) = 1,\, 2 \text{ and } 4 \text{, resp.}
  \end{equation*}
  The corresponding Reidemeister number is given by
  \begin{equation*}
    R_N = 2, 1 \text{ and } 1 \text{, resp.}	
  \end{equation*}
  The canonical fibration
  \begin{equation*}
    p : S^{n+d-1} \longrightarrow \KP(n')
  \end{equation*}
  will play a crucial r\^ole.
\end{exa}


\begin{thm}\label{thm:1.18}
  \newlength{\identlegth}
  \setlength{\identlegth}{\the\parindent}
  Assume that
  $ n' \geq 2 \text{ or } \pi_{m-1}(S^{d-1}) = 0 $.
  \\Then:\\
  (i) Given
  $ [f_i] \in \pi_m(\KP(n')) $,
  there exists a unique homotopy class
  $ [\widetilde{f}_i] \in \pi_m(S^{n+d-1}) $
  such that
  $ [f_i] - p_*[\widetilde{f}_i] $
  lies in the image of
  $ \pi_m(\KP(n') \backslash \{*\}),\, i= 1,2 $.
  (Since this image is isomorphic to
  $ \pi_{m-1}(S^{d-1}) $
  we may assume that
  $ \widetilde{f}_i $
  is a genuine lifting of
  $ f_i \text{ when } \K = \R \text{ or when } m > 2 \text{ and } \K = \C $).\\
  (ii) Given
  $ [f_1], [f_2] \in \pi_m(\KP(n')) $,
  assume that 
  $ (f_1, f_1) $
  is loose (cf. (\ref{def:1.4}); e.g this holds always when
  $ ( \K, m, n') $
  satisfies the assumptions of proposition \ref{prop:1.20} below).
  \\Then:
  \begin{xalignat*}{3}
    \titem&)\;& MCC(f_1, f_2) &= R_N \cdot MCC(\widetilde{f}_1, \widetilde{f}_2) 			&(&= 0 \longleftrightarrow \widetilde{f}_1 \sim \widetilde{f}_2). \\
    \titem&)\;& N^{\#}(f_1, f_2) &= R_N \cdot N^{\#}(\widetilde{f}_1, \widetilde{f}_2) 			&(&= 0 \longleftrightarrow \widetilde{f}_1 \sim \widetilde{f}_2). \\
    \titem&)\;& \widetilde{N}(f_1, f_2) &= R_N \cdot \widetilde{N}(\widetilde{f}_1, \widetilde{f}_2) 	&(&= 0 \longleftrightarrow \Gamma(\widetilde{f}_1) = \Gamma(\widetilde{f}_2)\,);
    \end{xalignat*}
    \leftskip=\identlegth here
    \[\Gamma \,\, {:=}\, \oplus E^{\infty} \scr \circ \dis \gamma_{k} \,:\, [S^m, S^{n+d-1}] \cong \pi_m(S^{n+d-1}) \longrightarrow \underset{k \geq 1}{\oplus} \pi_{m-1-k(n+d-2)}^S \]
    where $ E^{\infty} \scr \circ \dis \gamma_{k} $ denotes the stabilized $k^{\text{th}}$ Hopf--James invariant homomorphism.
    \begin{xalignat*}{3}
      \titem&)\;& N(f_1, f_2) &= R_N \cdot N(\,(E^{n-d}(h_{\K})) \scr \circ \dis \widetilde{f}_1,\, (E^{n-d}(h_{\K})) \scr \circ \dis \widetilde{f}_2\;) 	&& \\
      &&\; (&= 0 \;\longleftrightarrow\; E^{\infty}(h_{\K}) \cdot E^{\infty}([\widetilde{f}_1]-[\widetilde{f}_2]) \,=\, 0\; );	&&
    \end{xalignat*}
    \leftskip=\identlegth here
    $ h_{\K} \colon S^{2d-1} \longrightarrow \KP(1) = S^d $
    denotes the canonical projection ('Hopf map'); its infinite suspension
    $ E^{\infty}(h_{\K}) $
    represents
    $ 2 \in \pi_0^S = \Z $,
    the generators
    $ \eta \in \pi_1^S \cong \Z_2 \text{ or } \nu \in \pi_3^S \cong \Z_{24} $		
    according as
    $\; \K = \R, \C \text{ or } \H, \text{resp.} $
    \begin{xalignat*}{3}
      \titem&)\;& N^{\Z}(f_1, f_2) &= \begin{cases}
                                             R_N & \text{ if } m = n \text{ and } \widetilde{f}_1 \not\sim \widetilde{f}_2; \\
                                             0   & \text{ otherwise. }
                                      \end{cases} 		&&
    \end{xalignat*}
  
  In particular, the minimum number
  $ MCC(f_1, f_2) $
  and all four Nielsen numbers of
  $ (f_1, f_2) $
  can assume only the values
  $ 0 \text{ and } R_N $.
  Therefore these numbers are completely determined by the vanishing criteria spelled out above.
\end{thm}


If
$ n' = 1 \text{ then } \KP(n') $
is a sphere and these numbers are already known whether
$ \pi_{m-1}(S^{d-1}) $
vanishes or not (see our example \ref{exa:1.15}).	
In particular, we can deduce the following 'Wecken theorem'.

\begin{cor}\label{cor:1.19}
  If 
  $ n' $
  is odd and
  $ \K = \R \text{ or } \C $,
  or if
  $ \,n' \equiv 23 \, \mathrm{mod}\, 24 $
  and
  $ \K = \H $,
  then 
  \begin{equation*}
    MCC(f_1, f_2) = N^{\#}(f_1, f_2)
  \end{equation*}
  for all maps
  $ \, f_1, f_2 \, \colon \, S^m \longrightarrow \KP(n') \text{ where } m,n' \geq 1 $.
\end{cor}

A key ingredient in the proof of this corollary is the

\begin{prop}\label{prop:1.20}
  Given
  $\, m, n' \,\geq\, 1 \,$,
  assume that
  \begin{equation*}
    \K = \R \text{ or } \C,\, n' \equiv 1\, (2), 
    \quad \text{ or }\quad \K = \H,\, n' \equiv 23\, (24), 
    \quad \text{ or }\quad n = d n' \leq 3\,.
  \end{equation*}
  If
  $\, (m, n) \neq (2, 2) \,$
  then for all maps
  $ f: S^m \longrightarrow \KP(n') \,$
  the pair
  $ (f, f) $
  is loose.
  (In fact,
  $\, (f,f) $
  is even loose by small deformation, i.e. there exists an arbitrarily close approximation
  $ f' \text{ of } f $
  such that the pair
  $ (f,f') $
  is coincidence free).
\end{prop}

\begin{rem_cor}\label{rem:1.21}
  The assumptions in this proposition cannot be dropped.
  Indeed, consider the fiber projection
  $ p \,\colon\, S^{d(n'+1)-1} \longrightarrow \KP(n') $.
  If
  $ \K = \R \text{ or } \C $
  and
  $ n' $
  is even, or if
  $ \K = \H $
  and
  $ n' \not\equiv 23 \mod 24 $,
  then the pair
  $ (p, p) $
  is not loose.
  The somewhat unexpected claim for the quaternions is due to their noncommutativity on the one hand, and to the order of the stable 3--stem
  $ \pi_3^S \cong \Z_{24} $
  on the other hand.
  (In [K6],	
  Proposition 1.17 
  and the last three lines in Example 4.4 
  have to be corrected accordingly when
  $ \K = \H $).
\end{rem_cor}

When
$ \K = \R \text{ or } \C \,$
proposition \ref{prop:1.20} holds due to the fact that
$ \K^{n'+1} $
allows multiplication with the element
$ (0, 1) $
of the division algebra
$ \K \times \K $
of complex or quaternionic numbers, resp;
we can use the resulting tangential vector field on the unit sphere
$ S^{n+d-1} $
to push each fiber of
$ p_{\K} $
away from itself.
\hfill$\Box$

\begin{rem}\label{rem:1.22}
  The claims I)--V) in theorem \ref{thm:1.18} still hold for even
  $ n' $
  since we assume that
  $ (f_1, f_1) $
  is loose (and consequently
  $ (\widetilde{f}_1, \widetilde{f}_1) $
  is also loose and hence
  $ \widetilde{f}_1 \sim a \comp \widetilde{f}_1 $,
  where
  $ a $
  denotes the antipodal map).
  However, when
  $ n' $ is even this assumption often fails to hold---sometimes with striking consequences.
\end{rem}

Here we mention only one of many such cases:

\begin{exa}\label{exa:1.23} {\bfseries $ \mathbf{n = 16, 32} \text{ or } \mathbf{64,\, m = 2n-2, \,\mathbf{\mathbb{K}} = \mathbf{\mathbb{R}}.} $}	
  In these three dimension settings (and possibly also when
  $ n = 128 \text{ and } m = 254 $) 
  there exists a map
  $ f \colon S^m \longrightarrow \RP(n) $
  such that
  \[
    2 \,=\, R_N \,\neq\, 1 \,=\, MCC(f, f) \,\neq\, N^{\#}(f, f) \,=\, MCC(\widetilde{f}, \widetilde{f}) \,=\, 0	
  \]
  (cf. [K6], 1.27 or [KR], 1.13).	
  In particular, corollary \ref{cor:1.19} and several central claims in theorem \ref{thm:1.18}(ii) fail to hold.
\end{exa}

This 'non--Wecken' result is due to the existence of Kervaire invariant one elements in
$ \pi_m(S^n) $.
Their important r\^ole in coincidence theory was first pointed out in [GR2]	
and studied systematically in [K6]	
and [KR].	
In fact, [KR]	
discusses also coincidences of maps into arbitrary spherical space forms
$ N = S^n/G $
(i.e. orbit manifolds of free smooth actions of any finite group
$ G $
on
$ S^n $)
very carefully.
\hfill$\Box$

\begin{question}
  Is
  $\, \MCC \,\equiv\, N^{\#} \,$
  whenever
  $\, \K \,=\, \C \,$
  or
  $\, \H $?
\end{question}

As we have seen non--Wecken results of the form
$\, MCC \,\not\equiv\, N^{\#} \,$
can occur only when
$\, \K \,=\, \R \,$
or
$\, \C \,$
and
$ n' $
is even, or when
$ \K = \H \text{ and } n' \not\equiv 23 \mod 24 $.
In contrast, pairwise differences between our four types of Nielsen numbers are very common (this is already indicated in \ref{equ:1.14}) and lead to non-Wecken theorems of the form
$ \MCC \not\equiv \widetilde{N} \text{ or } \MCC \not\equiv N $.
(In fact, we do not expect interesting Wecken theorems
$ \MCC \equiv N^{\Z} $
at all in higher codimensions
$ m-n > 0 $).
Thus the following consequence of our discussion and, in particular, of theorem \ref{thm:1.18} underlines the importance of the Nielsen number
$ N^{\#} \,$
(based on {\itshape nonstabilized normal bordism} theory) when we try to compute minimum numbers.

\begin{cor}\label{cor:1.24}
  Let
  $ \K $
  be the field 
  $ \R,\, \C \text{ or } \H $
  (with real dimension
  $ d=1,\, 2 \text{ and } 4 $,
  resp.) and let
  $ n' $
  be an (even or odd) integer such that
  $ n \;{:=}\; dn' \geq 2 $.
  Then
  \begin{enumerate}[label=\alph*.)]
    \item
      $ N^{\#} \not\equiv \widetilde{N} $
      except possibly when
      $ \K = \R $
      and
      $ n \geq 12 $
      is even;
    \item
      $ \widetilde{N} \not\equiv N $;
      and
    \item
      $ N \not\equiv N^{\Z} \,$
      except precisely when
      $\, \K \,=\, \R \,$
      and
      $\, n \,=\, 2 $.
  \end{enumerate}
  Here
  $ N^{\#} \not\equiv \widetilde{N} $
  means that there exists
  $ m \in \Z $
  and maps
  $ f_1, f_2 \colon S^m \longrightarrow \KP(n') $
  such that
  $ N^{\#}(f_1, f_2) \neq \widetilde{N}(f_1, f_2) $,
  and similarly for the claims
  $ \widetilde{N} \not\equiv N $
  and
  $ N \not\equiv N^{\Z} $
  (possibly with different choices of
  $ m$).
  
  However, for
  $\, n \,=\, 1 \,$
  and all
  $\, m \geq 1 \,$
  the minimum numbers
  MC and MCC and all four Nielsen numbers agree for arbitrary pairs of maps
  $\, f_1, f_2 \,\colon\, S^m \longrightarrow \KP(n') $.
\end{cor}

To get a more precise picture we may want to fix not only
$ \K $
and
$ n' $
but also
$ m $,
and ask whether e.g.
$ N^{\#} \,\equiv\, \widetilde{N} $
(or
$ N^{\#} \,\not\equiv\, \widetilde{N} $)
in this context, i.e. whether (or not)
$ N^{\#}(f_1, f_2) \,=\, \widetilde{N}(f_1, f_2) $
for all maps
$ f_1, f_2 \,\colon\, S^m \longrightarrow \KP(n') $.
For this and similar comparisons involving also the Nielsen numbers
$ N $
and
$ N^{\Z} $
consider the commuting diagram of homomorphisms
\stepcounter{thm}
\addtocounter{equation}{8}
\begin{equation}
  \begin{aligned}\label{equ:1.25}
    \xymatrix {
      & \underset{k \geq 1}{\oplus} \pi_{m-1-k(n+d-2)}^S \ar[dd]^{E^{\infty}(h_{\K}) \cdot \text{ first projection} } \\
      \pi_m(S^{n+d-1}) = [S^m, S^{n+d-1}] \ar[ru]^{\Gamma} \ar[rd]_{E^{\infty}(h_{\K}) \cdot E^{\infty}\;\;\; } \\
      & \pi_{m-n}^S
    }
  \end{aligned}
\end{equation}
(where
$ \Gamma $
and
$\, E^{\infty}(h_{\K}) \,=\, 2,\, \eta $
or
$ \nu \,\in\, \pi_{*}^S $
are described in \ref{thm:1.18}(ii), III and IV).\\
We have
\stepcounter{thm}
\begin{equation}\label{equ:1.26}
  \{0\} \quad \underset{(a)}{\subseteqq} \quad \Ker \Gamma \quad \underset{(b)}{\subseteqq} \quad \Ker(E^{\infty}(h_{\K}) \cdot E^{\infty}) \quad \underset{(c)}{\subseteqq} \quad \pi_m(S^{n+d-1})
\end{equation}

Now assume that
\begin{enumerate}[label=(\roman*)]
  \item 
    $\, n' \geq 2 \;$
    or
    $\; \pi_{m-1}(S^{d-1}) = 0; \qquad \qquad $
    and
  \item for all maps
    $\, f \,\colon\, S^m \longrightarrow \KP(n') \,$
    the pair
    $ (f, f) $
    is loose.
\end{enumerate}
Then (according to theorem \ref{thm:1.18})
$ \; N^{\#} \; \equiv \;\widetilde{N} \;$
(or
$ \; \widetilde{N} \; \equiv \; N,\; $
or
$ \; N \; \equiv \; 0,\; $
resp.) if and only if we have a full equality---and not just an inclusion---at (a) (or (b), or (c), resp.) in diagram \eqref{equ:1.26};
when
$ m \neq n $
an equality at (c) is also equivalent to
$ N \,\equiv\, N^{\Z} $.

\begin{exa}\label{exa:1.27}
  {\bfseries n = 2.}
  When
  $\, \K = \C,\, n' = 1 $,
  then
  $ \KP(n') \cong S^2 $
  and hence
  $ \MCC \,\equiv\, N^{\#} $
  (cf. \ref{exa:1.15}).
  It is not hard to compare the Nielsen numbers for low values of
  $ m $
  (using standard techniques of homotopy theory such as EHP--sequences, and the tables of Toda [T]):	
  \begin{xalignat*}{10}
    m &= 2:	&	&N^{\#}	&	&\equiv		& &\widetilde{N}	&	&\equiv		&	&N	&	&\equiv		&	&N^{\Z}	&	&\not\equiv	&	&0 \\
    m &= 3:	&	&N^{\#}	&	&\equiv		& &\widetilde{N}	&	&\not\equiv	&	&N	&	&\not\equiv	&	&N^{\Z}	&	&\equiv		&	&0 \\
    m &= 4,5:	&	&N^{\#}	&	&\equiv		& &\widetilde{N}	&	&\equiv		&	&N	&	&\not\equiv	&	&N^{\Z}	&	&\equiv		&	&0 \\
    m &= 6,7,8:	&	&N^{\#}	&	&\equiv		& &\widetilde{N}	&	&\not\equiv	&	&N	&	&\equiv		&	&N^{\Z}	&	&\equiv		&	&0 \\
    m &= 9:	&	&N^{\#}	&	&\not\equiv	& &\widetilde{N}	&	&\equiv		&	&N	&	&\equiv		&	&N^{\Z}	&	&\equiv		&	&0
  \end{xalignat*}
  Here we get e.g. a Wecken theorem of the form
  $ \MCC \,\equiv\, N $
  precisely when
  $ m=2,\, 4 $
  or
  $ 5 $,
  and no Wecken theorem of the form
  $ \MCC \,\equiv\, \widetilde{N} $
  when
  $ m = 9 $.
  \hfill$\Box$
\end{exa}

When
$\, \K \,=\, \R \,$
and
$\, n' \,=\, 2 \,$
and we consider maps
$\, f_1, f_2 \,\colon\, S^m \,\longrightarrow\, \RP(2) $,
we can compare the Nielsen numbers
$\, N^{\#}, \, \widetilde{N}, \, N \,$
and
$\, N^{\Z} \,$
of
$\, (f_1, f_2) \,$
between themselves, but also with the Nielsen numbers of the liftings
$\, \widetilde{f}_1, \widetilde{f}_2 \,\colon\, S^m \,\longrightarrow\, S^2 \cong \CP(1) $.
It follows from theorem \ref{thm:1.18}(ii) that
$\, N^{\#}(f_1, f_2) \,\equiv\, 2 \cdot N^{\#}(\widetilde{f}_1, \widetilde{f}_2) \,$
and
$\, \widetilde{N}(f_1, f_2) \,\equiv\, 2 \cdot \widetilde{N}(\widetilde{f}_1, \widetilde{f}_2) $,
but
$\, N(f_1, f_2) \,\equiv\, N^{\Z}(f_1, f_2) $,
and, if
$\, m \geq 3 \,$
and
$\, E^{\infty}(\pi_m(S^2)\,) \,\neq\, \{0\} $,
then
\begin{equation*}
  N(\widetilde{f}_1, \widetilde{f}_2) \,\not\equiv\, 2 \cdot N(f_1, f_2) \,\equiv\, 0
\end{equation*}
\hfill$\Box$

\vspace{2ex}
For more background and some of the many further aspects of Nielsen fixed point and coincidence theory or normal bordism techniques consult e.g. also the papers
[B2], [BS], [BGZ], [C], [D], [DG], [GR1], [HQ], [K1], [K4], [K5] and [S]	
listed in our references (no claim to completeness!).

\section{Coincidences in projective spaces}

In this section we prove theorem \ref{thm:1.18}.

  (i) Note that
  $ \KP(n'-1) $
  is a deformation retract of the punctured projective space
  $ \KP(n') - \{*\}$.
  Moreover the fiber map
  $ p:S^{n+d-1} \rightarrow \KP(n') $
  with nulhomotopic fiber inclusion i gives rise to the commuting diagram
  \begin{equation}
    \begin{aligned}\label{equ:2.1}
      \xymatrix  @R=0.4cm @C=0.4cm  {
        \pi_m( S^{n-1} )	\ar[r]^-{p\vert_*} \ar[d]^0	&
        \pi_m( \KP(n'-1) )	\ar[r]^-{\partial\vert} \ar[d]^-{\incl_*}	&
        \pi_{m-1}( S^{d-1} )	\ar[r]^-{i\vert_*} \ar@{-->}[ld] \ar@{=}[d]	&
        \pi_{m-1}( S^{n-1} )	\ar[d]	\\
        \pi_m( S^{n+d-1} )	\ar@{^{(}->}[r]^-{p_*}	&
        \pi_m( \KP(n') )	\ar@{->>}[r]	&
        \pi_{m-1}( S^{d-1} )	\ar[r]^-0	&
        \pi_{m-1}( S^{n+d-1} )
      }
    \end{aligned}
  \end{equation}
  Thus
  $ p_* $
  is injective. If
  $ n' \geq 2 $,
  then the fiber inclusion
  $ i\vert : S^{d-1} \subset S^{n' d -1} $
  is also nulhomotopic and
  $ \incl_* $
  factors through the boundary epimorphism
  $ \partial\vert $
  ; this yields a splitting of the lower horizontal sequence and an isomorphism
  $ \incl_*(\pi_m (\KP(n')-\{*\}) ) \cong \pi_{m-1}( S^{d-1} ) $.
  If
  $ n' = 1 $
  and
  $ \pi_{m-1}( S^{d-1} ) $
  vanishes, then so do the image of
  $ \incl_* $
  and the cokernel of
  $ p_* $.\\
  
  (ii) Choose basepoints
  $ *_1 \neq *_2 $
  in
  $ N= \KP(n'),\, \K = \R,\, \C \text{ or } \H $.
  Given classes
  $ [f_i], [l_i] \in \pi_m( N, *_i), i=1,2 $,
  such that
  $ ( l_1, l_2 ) $
  is loose (in the basepoint free sense), it is easy to see that the pairs
  $ ( [f_1],[f_2] ) $
  and
  $ ( [f_1] + [l_1],[f_2]+[l_2] ) $
  have the same minimum and Nielsen numbers (cf. also the appendix in [K3]).	
  E.g. if
  $ [f_i] = [ p \comp \widetilde{f}_i ] - [ l_i ] $,
  where
  $ [l_i] \in \incl_*( \pi_m( \KP(n') - *_{i \pm 1} )), i=1,2 $,
  as in claim (i) 
  of our theorem, then
  \[
    ( [l_1],[l_2] ) = ( [l_1],[*_2] ) + ( [*_1],[l_2] )
  \]
  is loose. Therefore we may assume henceforth in our proof that
  $ [f_i] = [ p \comp \widetilde{f}_i ],\, i=1, 2 $.

  Next choose
  $ [f_1'] \in \pi_m(\KP(n'),*_2 ) $
  such that
  $ f_1' \sim f_1 $
  (just the basepoint behaviour is modified, e.g. by an isotopy of
  $ N $).
  Then
  $ ( l_1, l_2 ) {:=} ( f_1, f_1' ) $
  is loose by assumption.
  We define
  \begin{equation} \label{equ:2.2}
    [f]\, {:=}\, [f_2] - [f_1']
  \end{equation}
  and we see (as above) that the pairs
  $ ( [f_1],[f_2] ) \text{ and } ( ( [f_1]-[f_1],[f_2]-[f_1']) \, = \, (* \,{:=}\, *_1,\, [f] ) $
  have the same minimum and Nielsen numbers.
  Thus we need to consider only pairs of the form
  $ (\,*\,=p(\widetilde{x}) ,\, f = p \comp \widetilde{f} ) $
  in our proof.
  
  Then the vanishing part of claim I in theorem \ref{thm:1.18} is obvious: if
  $ \MCC(*, f) = 0 $
  and hence
  $ f $
  can be deformed into
  $ \KP(n') \backslash \{*\} $
  then the lifting
  $ \widetilde{f} $
  is homotopic to a map into
  \begin{equation*}
    S^{n+d-1} \backslash p^{-1}(\{*\}) \; \subset \; S^{n+d-1} \backslash\{\widetilde{*}\} \;\sim\; \{ - \widetilde{*} \};
  \end{equation*}
  in turn, if
  $ \MCC( \widetilde{*}, \widetilde{f}) = 0 \,$
  or, equivalently,
  $ \widetilde{f} $
  is nulhomotopic then so is 
  $ f = p \scr \circ \dis \widetilde{f} $
  and
  $ MCC(*, f) = 0 $.
  
  Moreover let us recall that the pairs
  $ (*, f) $
  and
  $ (f, *) $
  have equal Nielsen numbers
  $ N^{\#} $
  and
  $ \widetilde{N} $
  (cf. [K3], 1.2(ii)).	
  However,
  $ N(*, f) $
  may differ from
  $ N(f, *) $
  but is easier to describe (due to our framing convention in the construction of
  $ \omega$--invariants,
  cf. (\ref{equ:1.6},ii)).
  
  Let us compare
  $ \widetilde{N}(f,*) \text{ to } \widetilde{N}(\widetilde{f},\widetilde{*}) $.
  (The corresponding discussion of
  $ N^{\#}(f,*) $
  vs.
  $ N^{\#}(\widetilde{f},\widetilde{*}) $
  was carried out in greater generality in the proof of theorem 6.5 in [K3]). 
  Consider the diagram of homomorphisms
  \begin{equation} \label{equ:2.3}
    \begin{aligned}
      \xymatrix{
        **[r] \pi_m( S^{n+d-1} ) \ar[d]^{\widetilde{\deg}_Q} \ar[r]^{p_*} &
        **[l] \pi_m(\KP(n')) \ar[d]^{\widetilde{\deg}_N} \\
       **[r] \pi^S_{m-n-(d-1)}(\Omega( Q,\widetilde{y_0} )) \ar@<-2pt>[r]_-{\beta} &
       **[l] \pi^S_{m-n}(\Omega( N,y_0 )) \ar@<-2pt>[l]_-{\alpha}
      }
    \end{aligned}
  \end{equation}
  where we write
  $ Q {:=} S^{n+d-1} $ and 
  $ N {:=} \KP(n') $
  for brevity,
  $ p_* $ is induced by the fiber projection,
  $ y_0 \, {:=} \, p(\widetilde{y_0}) $,
  and the vertical arrows are defined by (\ref{def:1.11}).
  We will now describe the homomorphisms
  $ \alpha \text{ and } \beta $.
  
  Given an element
  $ c \in \pi_{m-n}^S(\Omega(N,y_0)) $,
  interpret it---via the Pontrjagin--Thom procedure---as a framed bordism class of a framed (= stably parallelized)
  $ (m-n)$--dimensional
  manifold
  $ C $, equipped with a map
  $ \widetilde{g} \colon C \rightarrow \Omega N $.
  The corresponding evaluation map
  $ C \times I \rightarrow N $
  lifts to a homotopy
  $ \widetilde{G} $ from the constant map at the point
  $ \widetilde{y_0} \in Q $
  to a map
  $ \widetilde{G}_1 \colon C \rightarrow F $
  into the fiber
  $ F = p^{-1}(\{y_0\}) \subset Q $.
  We may assume
  $ \widetilde{G}_1 $
  to be smooth, with regular value
  $ \widetilde{y_0} $.
  Equip
  $\, C' \,{:=}\; \widetilde{G}_1^{-1}(\{y_0\}) $
  with the map
  $ \widetilde{g}' \colon C' \rightarrow \Omega( Q,\widetilde{y_0} ) $
  which corresponds to
  $ \widetilde{G}\vert C' \times I $.
  Moreover compose the natural trivialization of the normal bundle of
  $ C' \text{ in } C $
  (given by the tangent map of
  $ \widetilde{G}_1 $)
  with the automorphism of
  $ C' \times T_{\widetilde{y_0}} F $
  which is determined by the homotopy
  $ \widetilde{G} \vert C' \times I $
  and the tangent bundle along the fibers of
  $ p $
  (cf. [K2], 3.1). 
  The resulting framed bordism class
  $ [C', \widetilde{g}'] $ defines $ \alpha(c) $
  (again via Pontrjagin--Thom).
  We have
  \begin{equation} \label{equ:2.4}
    \widetilde{\deg}_Q \;=\; \pm \alpha \scr \circ \dis \widetilde{\deg}_N \scr \circ \dis  p_*
  \end{equation}
  since
  $ \widetilde{\deg}_N \scr \circ \dis p_* $ and
  $ \widetilde{\deg}_Q $
  correspond to taking the inverse image of a fiber and of a point in
  $ Q $,
  resp.
  
  Since
  $ n' \geq 1 $
  there exists a homotopy 
  \begin{equation*}
    \widetilde{J} \colon ( F, \widetilde{y}_0 ) \times I \rightarrow ( Q, \widetilde{y}_0 )
  \end{equation*}
  from the constant map at
  $ \widetilde{y}_0 $
  to the inclusion of the fiber
  $ F \cong S^{d-1} $ into
  $ Q = S^{dn'+d-1} $.
  
  Given an element
  $ c' \in \pi^S_{m-n-(d-1)}( \Omega(Q,\widetilde{y}_0) ) $,
  describe it by a framed manifold
  $C'$,
  together with a map
  $ \widetilde{g}' \colon C' \rightarrow \Omega(Q,\widetilde{y}_0) $.
  Endow
  $ F \cong S^{d-1} $
  with the (left invariant) Lie group framing and
  $ C \,{:=}\, C' \times F $
  with the resulting product framing. Moreover let
  \begin{equation*}
    \widetilde{g} \;\colon\; C = C' \times F \longrightarrow \Omega(N,y_0)
  \end{equation*}
  be given by the loops in
  $ N $
  which concatenate
  $ p \scr \circ \dis \widetilde{g}' $
  with the adjoint of
  $ p \scr \circ \dis \widetilde{J} $.
  We obtain
  $ \beta(c') $
  by applying the Pontrjagin--Thom isomorphism to the framed bordism class of
  $ (C, \widetilde{g} ) $.
  
  This definition of
  $ \beta $
  mimics the transition from
  $ \widetilde{\deg}_Q $ to
  $ \widetilde{\deg}_N $
  where the inverse image of a point
  $ \widetilde{*} \in Q $
  is replaced by the inverse image of the whole fiber containing
  $ \widetilde{*} $.
  We obtain
  \begin{equation}\label{equ:2.5}
    \beta \scr \circ \dis \widetilde{\deg}_Q \;=\; \pm \widetilde{\deg}_N \scr \circ \dis p_*
  \end{equation}
  and
  $ \alpha \scr \circ \dis \beta = \pm \id $.
  Therefore
  $ \beta $
  is injective.
  
  It follows from (\ref{equ:2.2}), (\ref{equ:2.4}) and (\ref{equ:2.5}) that
  $ \widetilde{N}(f_1, f_2) = \widetilde{N}(f, *) $
  (and hence
  $ \widetilde{\deg}_N(\{f\}) \,$)
  vanishes if and only if
  $ \widetilde{\deg}_Q([\widetilde{f}]) $
  does or, equivalently,
  $ \Gamma([\widetilde{f}]) = 0 $
  (cf. theorem 1.14 in [K2]). 
  This establishes the vanishing criterion in claim III 
  of our theorem (and similarly in claim II 
  since
  $ \deg^{\#}_{Q} $
  is injective).
  
  The looseness obstructions
  $ \omega(f_1, f_2) $ and
  $ \omega(\widetilde{f}_1,\widetilde{f}_2) $
  are obtained from
  $ \widetilde{\omega}(f_1, f_2) $ and
  $ \widetilde{\omega}(\widetilde{f}_1,\widetilde{f}_2) $,
  resp., by forgetting the maps into loopspaces.
  According to the framing convention emboddied in the construction of our
  $ \omega$--invariants (cf. e.g. (\ref{equ:1.6},ii) or [K3], 
  formulas (6) and (22)),
  $ \,\omega(\widetilde{f}_1, \widetilde{f}_2) \,=\, \omega(\widetilde{*},\widetilde{f}) $
  corresponds---via the Pontrjagin--Thom isomorphism---to  the bordism class of the (generic) inverse image manifold
  $ \widetilde{f}^{-1}(\{\widetilde{*}\}) $,
  framed in the obvious fashion; thus
  \begin{equation*}
    \omega(\widetilde{f}_1, \widetilde{f}_2) \,=\, \omega(\widetilde{*},\widetilde{f}) \,=\, E^{\infty}([\widetilde{f}]) \;\in\; \pi^S_{m-n-(d-1)}.
  \end{equation*}
  Similarly
  $ \omega(*,h_{\K}) =E^{\infty}(h_{\K}) \,\in\, \pi^S_{d-1} $
  corresponds to the framed bordism class of the fiber of the Hopf map
  $ h_{\K} $,
  i.e. to
  $ [S^{d-1}] $
  where the Lie group
  $ S^{d-1} $
  is endowed with its left invariant framing.
  On the other hand it follows as in the previous discussion that
  $ \omega(*, f) $
  corresponds to the product of
  $ \widetilde{f}^{-1}(\{\widetilde{*}\}) $
  with
  $ S^{d-1} $.
  Therefore
  \begin{equation*}
    \pm\, \omega(f_1, f_2) \,=\, E^{\infty}([h_{\K}]) \cdot E^{\infty}[\widetilde{f}] \,=\, E^{\infty}( E^{n-d}([h_{\K}]) \comp [\widetilde{f}]).
  \end{equation*}
  
  In contrast to the situation in the cases II and III,
  $ \omega(*, f) $ and
  $ \omega(\widetilde{*}. \widetilde{f}) $
  need not be equally strong since here we lack maps into the loop space
  $ \Omega N $
  and hence a homomorphism
  $ \alpha $ as in (\ref{equ:2.4}).
  This explains the different form of claim IV.
  
  Finally recall that the values of MCC and the Nielsen numbers are bounded from above by the Reidemeister number
  $ R_N $ (cf. \ref{equ:1.12} and \ref{equ:1.14});
  we know this even when
  $ n=2 $ since it is true for spheres (use surgery; cf. also \ref{exa:1.15}).
  In fact, 
  $ R_N $
  is the only possible nontrivial value. This holds obviously when
  $ N $
  is a sphere or a complex or quaternionic projective space since
  $ R_N = 1 $
  in this case. If
  $ N = \RP(n) $
  and 
  $ \widetilde{*} $
  is a regular value of a smooth map
  $ \widetilde{f} : S^m \longrightarrow S^n $ and
  $ * \;{:=}\; p(\widetilde{*}) = p(- \widetilde{*}) \in \RP(n) $,
  then the coincidence manifold
  $ C(*, f) = f^{-1}(\{*\}) $
  consists of the two Nielsen classes
  $ \widetilde{f}^{-1}(\{\widetilde{*}\}) $ and
  $ \widetilde{f}^{-1}(\{-\widetilde{*}\}) $
  which may be assumed to be connected and which---for each of the Nielsen numbers
  $ N^{\#},\, \widetilde{N},\, N \text{ and } N^{\Z} $---are
  simultaneously either essential (or not) according as
  $ \widetilde{f}(\{\widetilde{*}\}) $
  contributes nontrivially (or not) to the Nielsen number in question for the lifted pair
  $(\widetilde{*}, \widetilde{f}) $.
  E.g. if 
  $ \MCC(*, f) \neq 0 $,
  then (by the vanishing criterion in case I)
  $\; \MCC(\widetilde{*}, \widetilde{f}) = N^{\#}(\widetilde{*}, \widetilde{f}) $
  (cf. \ref{exa:1.15})
  is nontrivial and both Nielsen classes of
  $ (*, f) $
  contribute nontrivially to
  $\, N^{\#}(*, f) \leq \MCC(*, f) \leq R_N = 2 $;
  thus 
  $ \MCC(*, f) = R_N $.
  
  The full claims I--IV in theorem \ref{thm:1.18} follow now from the vanishing criteria, and so does claim V.
  Indeed, since
  $\, m \geq 2 \,$
  and
  $\, n \neq 0 $, 
  $ N^{\Z}(f_1, f_2) \,$
  and
  $\, \omega^{\Z}(f_1, f_2) \,\in\, H_{m-n}(S^m; \Z) \,$
  can be nontrivial only when
  $ m = n $
  and therefore
  \begin{equation*}
    N^{\Z}(f_1, f_2) \,=\, N^{\#}(f_1, f_2) \,=\, R_N \cdot N^{\#}(\widetilde{f}_1, \widetilde{f}_2) \,\neq\, 0
  \end{equation*}
  and
  $ \widetilde{f}_1 \,\not\sim\, a \cdot \widetilde{f}_2 $
  (cf. \ref{exa:1.15}).
  \hfill$\Box$

\section{Selfcoincidences}

In this section we prove proposition \ref{prop:1.20} and the claims in remark \ref{rem:1.21}.

If
$\, \K = \C \,$
and
$\, (m, n') = ( 2, 1) \,$,
then
$\, \KP(n') \cong S^2 \,$
and the pair
$\, (f, f) \;$
(where
$\, f \,\colon\, S^2 \longrightarrow S^2 \,$)
is loose if and only if 
$\, f \sim (\text{antipodal map }) \comp f $,
i.e.
$\, f \sim \text{constant map} $.
Thus we have to exclude the case
$\, (\K,\, m,\, n') \,=\, ( \C, 2, 1) \,$
from further discussions.

If
$\, m \,\neq\, n = 2 \,$,
then
$\, \pi_{m-1}(S^1) \,=\, 0\,$
and every map
$\, f \,$
from
$\, S^m \,$
to
$\, \CP(1) \,$
or
$\, \RP(2) \,$
lifts to the tangent circle bundle of this surface;
we can use the resulting vector field (parametrized by
$\, S^m \,$)
to 'push
$\, f \,$
away from itself'.

Now assume that
$ n' $
is odd,
$ \, \K = \R,\, \C \text{ or } \H $.
Then we can multiply the elements of
$ \K^{n'+1} $
on the left with the element
$ (0,1) $
of the division algebra
$ \K \times \K $
(of complex, quaternionic or octonic numbers, resp.).
Restriction to the unit sphere yields the selfmap
$ s \text{ of } S^{n+d-1} $,
described by
\begin{equation*}
  x = (x_1, x_2; x_3, x_4; \ldots ; x_{n'}, x_{n'+1}) \longrightarrow (-\bar{x}_2, \bar{x}_1; -\bar{x}_4, \bar{x}_3; \ldots) \in S^{n+d-1} \subset \K^{n'+1}
\end{equation*}
which is homotopic to the identity map 
$( S^{2d-1} $
being connected).
If
$\, s(x) = \lambda x \,$
for some
$\, \lambda \in \K $,
then
\[
  -\bar{x}_{2i} \,=\, \lambda x_{2i-1},\quad \bar{x}_{2i-1} \,=\, \lambda x_{2i},\quad x_{2i-1} \,=\, \bar{x}_{2i} \cdot \bar{\lambda}
\]
and hence
\begin{equation}\label{equ:3.1}
  \bar{x}_{2i} \,=\, - \lambda \bar{x}_{2i} \bar{\lambda} \qquad\text{for } i= 1,\ldots, (n'+1)/2.
\end{equation}

In case
$\,\K\,$
is commutative we conclude that
$\, ( 1+ \vert \lambda \vert^2 ) \vert x_{2i} \vert = 0 \,$
and hence
$\, x = 0 \,$
contradicting our assumption that
$\, x \,\in\, S^{n+d-1} $;
therefore
$\, s(x) \notin \K \cdot x \,$,
and
$ s $
gives rise to a nowhere vanishing vectorfield
$ v $
in the pullback
$\, p^{*}(T \,\KP(n')) \,$
over 
$\, S^{n+d-1} \,$
along which we can push each fiber of
$ p $
away from itself; thus
$\, f = p \comp \widetilde{f} \,$
(cf. \ref{thm:1.18} (i) ) has no coincidence with
$\, p \comp s \comp \widetilde{f} \,$
and
$\, (f, f) \,$
is loose.

This whole argument depends on
$\, \bar{x}_{2i} \cdot \bar{\lambda} \,$
being equal to
$\, \bar{\lambda} \cdot \bar{x}_{2i} \,$
in \eqref{equ:3.1}.
In 
$ \H $
this need not hold and
$ s(x) $
may lie in the line
$ \K \cdot x $
(e.g.
$ s(x) = ix $
when
$ x = (j,k;0,\ldots, 0 )$\,);
thus
$\, f = p \comp \widetilde{f} \,$
may have coincidences with
$\, p \comp s \comp \widetilde{f} \,$.
However, if
$\, n' \equiv 23\,(24) \,$
then it follows from formula (5.9) in [Ja], p. 38,	
that there still exists a selfmap
$\, s' \text{ of } S^{4(n'+1)-1} \subset \H^{n'+1} \,$
such that
$\, s'(x) \,\notin\, \H x \,$
for all
$\, S^{4(n'+1)-1} \,$;
the pair
$\, (p \comp s' \comp \widetilde{f},\, \widetilde{f}) \,$
of homotopic maps is coincidence free.
This establishes proposition \ref{prop:1.20}.
\hfill$\Box$\\

Whether
$\, \K = \R,\, \C \text{ or } \H $,
the following conditions are equivalent for all
$ [f] \,\in\, \KP(n'),\, n' \geq 2 $,
(cf. [K6], theorem 1.22):	
\begin{enumerate}
  \item
  	$ (f, f) $
  	is loose by small deformation;
  \item
  	$ (f, f) $
  	is loose (by any deformation);
  \item
  	$ f $
  	is not coincidence producing (i.e. there exists some map
  	$\, f' \,\colon\, S^m \rightarrow \KP(n') \,$
  	such that the pair
  	$\, (f, f') \,$
  	is loose; compare [BS]).	
\end{enumerate}

However these three conditions need not be equivalent for other target manifolds, not even for maps between spheres (cf. e.g. [GR1], [GR2] or [K6], corollaries 1.21, 1.28 and 1.30).

The claim in remark \ref{rem:1.21} follows from [Ja], formulas 5.8 and 5.9	
(compare also [DG], theorems 3.5 and 3.9).	
Indeed,
$\, (p, p) \,$
is loose if and only if the canonical fibration of the Stiefel manifold
$\, V_{n'+1,2}(\K) \,$
(of orthonormal $ 2$--frames in
$\, \K^{n'+1} \,$)
over the sphere
$\, S(\K^{n'+1}) \,=\, S^{d(n'+1)-1} \,$
allows a section.
\hfill$\Box$

\newcommand{\doubletilde}[1]{\overset{\approx}{#1}}
\section{Examples}

In this section we discuss example \ref{exa:1.27} and corollary \ref{cor:1.24}.

First consider the case
$\, n= 2 \,$.
Fix
$\, m \geq 1 \,$
and the surface
$\, \CP(1) \,$
or
$\, \RP(2) \,$
and compare the Nielsen numbers of all pairs of maps from
$\, S^m \,$
into this surface.
When
$\, m = 1 \,$
or
$\, 2 \,$
the four Nielsen numbers agree.

Thus assume that
$\, m > n = 2 \,$.
Then
$\, \pi_{m-1}(S^{d-1}) = 0 \,$.
Also
$\, N^{\Z} \equiv 0 \,$,
but the other Nielsen numbers may show interesting differences between each other or between the two target surfaces.
We can exploit the double r\^ole which
$\, \CP(1) = S^2 \,$
plays both as a base and as a total space of a canonical projection (real or complex 'Hopf map'); these induce the isomorphisms
\begin{equation}\label{equ:4.1}
  \begin{aligned}
    \xymatrix @R=1pt {
      \pi_m(S^3)	\ar[r]^-{\cong}	&
      \pi_m(S^2)	\ar[r]^-{\cong}	&
      \pi_m(\RP(2)\,)	\\
      \rotatebox{90}{$\in$}	&
      \rotatebox{90}{$\in$}	&
      \rotatebox{90}{$\in$}	\\
      [\overset{\approx}{f}_i]	&
      [\widetilde{f}_i]	&
      [f_i].
    }
  \end{aligned}
\end{equation}
Proposition \ref{prop:1.20} and theorem \ref{thm:1.18} apply fully to all maps
\begin{equation*}
  f_1, f_2 \,\colon\; S^m \longrightarrow \RP(2)
\end{equation*}
and to their liftings
$\, \widetilde{f}_1, \widetilde{f}_2 \,$
into
$\, S^2 \,$
and
$\, \doubletilde{f}_1, \doubletilde{f}_2 \,$
into
$\, S^3 \,$
(cf. \ref{equ:4.1}).
We get the following conclusions from \ref{thm:1.18}(ii),II--IV:
\stepcounter{countertitem}
\\\titem)
$\qquad
  \xymatrix  @R=1pt {
    N^{\#}(f_1, f_2)	\ar@{}@<-4pt>[r]^{=} &
    2 \cdot N^{\#}(\widetilde{f}_1, \widetilde{f}_2)	\ar@{}@<-4pt>[r]^{=} &
    2 \cdot N^{\#}(\doubletilde{f}_1, \doubletilde{f}_2) \\
    (\;=\, 0\;	\ar@{<->}[r] &
    \widetilde{f}_1 \sim \widetilde{f}_2	\ar@{<->}[r] &
    \doubletilde{f}_1 \sim  \doubletilde{f}_2 \;).
  }
$
\\In particular,
\begin{equation*}
  N^{\#}(f_1, f_2) \,\not\equiv\, 0 \,\longleftrightarrow\, N^{\#}(\widetilde{f}_1, \widetilde{f}_2) \,\not\equiv\, 0 \,\longleftrightarrow\, \pi_m(S^2)(\cong \pi_m(S^3)\,) \,\neq\, 0;
\end{equation*}
this condition is satisfied e.g. for
$\, 3 \leq m \leq 21 \,$
(cf. Toda's tables in [T], p.186).	
\\\titem)
$\qquad
  \xymatrix  @R=1pt {
    \widetilde{N}(f_1, f_2)	\ar@{}@<-4pt>[r]^{=} &
    2 \cdot \widetilde{N}(\widetilde{f}_1, \widetilde{f}_2)	\ar@{}@<-4pt>[r]^{=} &
    2 \cdot \widetilde{N}(\doubletilde{f}_1, \doubletilde{f}_2) \\
    (\;=\, 0\;	\ar@{<->}[r] &
    \Gamma(\widetilde{f}_1) = \Gamma(\widetilde{f}_2)	\ar@{<->}[r] &
    \Gamma(\doubletilde{f}_1) = \Gamma(\doubletilde{f}_2)\;).
  }
$
\\In particular,
\begin{equation*}
  \widetilde{N}(f_1, f_2) \,\equiv\, 0 \,\longleftrightarrow\, \widetilde{N}(\widetilde{f}_1, \widetilde{f}_2) \,\equiv\, 0 \,\longleftrightarrow\, \Gamma(\pi_m(S^2)\,) \,=\, 0 \,\longleftrightarrow\,
  \Gamma(\pi_m(S^3)\,) \,=\, 0,
\end{equation*}
and we can choose whether we want to compute
$\, \Gamma \,$
on
$\, \pi_m(S^2) \,$
or, equivalently, on
$\, \pi_m(S^3) $.
E.g. when
$\, m = 9 \,$
we have the choice between the homomorphisms
\begin{subequations}
  \begin{align}
    \Gamma \;\colon\; \pi_9(S^2) \cong \Z_3 &\,\longrightarrow\, \Z_{240} \oplus \Z_2 \oplus \{0\} \oplus \{0\} \oplus \Z_{24} \oplus \Z_2 \oplus \Z_2 \oplus \Z \label{equ:4.2a}\\
    \intertext{and}
    \Gamma \;\colon\; \pi_9(S^3) \cong \Z_3 &\,\longrightarrow\, \Z_2 \oplus \{0\} \oplus \Z_2 \oplus \Z \label{equ:4.2b}
  \end{align}
\end{subequations}
(cf. \ref{exa:1.15} and \ref{thm:1.18}(ii), as well as [T], p. 186);	
without any calculation we see from \ref{equ:4.2b} that
$\, \Gamma \equiv 0 \,$
and hence
\begin{equation}\label{equ:4.3}
  N^{\#} \;\not\equiv\; \widetilde{N} \;\equiv\; 0
\end{equation}
for maps from
$\, S^9 \,$
to
$\, \CP(1) \,$
or
$\, \RP(1) $.
Thus the Nielsen number
$\, N^{\#} \,$
based on {\itshape nonstabilized} normal bordism theory is strictly more powerful than the Nielsen number
$\, \widetilde{N} \,$
which embodies only a standard bordism approach.
\vspace{6pt}
\\\titem)
$\qquad
  \begin{aligned}[t]
    N(f_1, f_2) = 0 &\longleftrightarrow\; 2 \cdot E^{\infty}([\widetilde{f}_1] - [\widetilde{f}_2]) \,=\, 0, \qquad \text{but}\\
    N(\widetilde{f}_1, \widetilde{f}_2) = 0 &\longleftrightarrow\; E^{\infty}([\widetilde{f}_1] - [\widetilde{f}_2]) \,=\, 0\;.
  \end{aligned}
$
\\Since
$\, \widetilde{f}_i \,$
is obtained by composing
$\, \doubletilde{f}_i \,$
with the (complex) Hopf map,
$\, i = 1,2 \,$,
we see that
\begin{equation*}
  E^{\infty}([\widetilde{f}_1] - [\widetilde{f}_2]) \,=\, \eta \cdot E^{\infty}([\doubletilde{f}_1] - [\doubletilde{f}_2])
\end{equation*}
vanishes when multiplied with
$\, 2 $;
thus
\begin{subequations}
  \begin{align}
    N(f_1, f_2) &\,\equiv\, N^{\Z}(f_1, f_2) \label{equ:4.4a}\\
    \intertext{for all
    $\, m \geq 1 $.
    On the other hand,}
    N(\widetilde{f}_1, \widetilde{f}_2) &\,\not\equiv\, N^{\Z}(\widetilde{f}_1, \widetilde{f}_2) \label{equ:4.4b}
  \end{align}
\end{subequations}
precisely when
$\, m \geq 3 \,$
and
$\; E^{\infty}(\pi_m(S^2)\,) \;=\; \eta \cdot E^{\infty}(\pi_m(S^3)\,) \;\neq\; 0 \,$,
e.g. when
$\, m = 3,\, 4 \,$
or
$\, 5 \,$,
but not when
$\, 6 \,\leq\, m \,\leq\, 9 \,$.\\

Let us take a closer look at the case
$\, (m, n) \,=\, (3, 2) $.
Here
\begin{equation*}
  \Gamma \;\colon\; \pi_3(S^2) \;\longrightarrow\; \pi_1^S \oplus \pi_0^S \;=\; \Z_2 \oplus \Z
\end{equation*}
is given by the suspension homomorphism and the Hopf invariant and hence injective; alternatively
\begin{equation*}
  \Gamma \,=\, E^{\infty} \;\colon\; \pi_3(S^3) \,\longrightarrow\, \pi_0^S
\end{equation*}
is just the suspension isomorphism.
Both ways we see that
$\, \Ker \Gamma \,=\, 0 $.
Thus
$\, N^{\#}(\widetilde{f}_1, \widetilde{f}_2) \,\equiv\, \widetilde{N}(\widetilde{f}_1, \widetilde{f}_2) \,$
vanishes if and only if
$\, \widetilde{f}_1, \widetilde{f}_2 \,\colon\, S^3 \longrightarrow S^2 \,$
have equal Hopf invariants.
If these Hopf invariants agree only\hspace{-7pt}
$ \mod{2} \,$
then
$\, N(\widetilde{f}_1, \widetilde{f}_2) \,=\, 0 $.
Taking into account also the corresponding maps
$\, f_1, f_2 \,$
into
$\, \RP(2) \,$
we conclude that
\begin{equation}\label{equ:4.5}
  N^{\#}(f_1, f_2) \,\equiv\,
  \widetilde{N}(f_1, f_2) \,\equiv\,
  2 \cdot \widetilde{N}(\widetilde{f}_1, \widetilde{f}_2) \,\not\equiv\,
  2 \cdot N(\widetilde{f}_1, \widetilde{f}_2) \,\not\equiv\,
  N(f_1, f_2) \,\equiv\, 0.
\end{equation}

When
$\, m \,=\, 4 \text{ or } 5 \,$
then
$\, E^{\infty},\, \eta \cdot E^{\infty} \,$
and
$\, \Gamma \,$
are injective on
$\, \pi_m(S^3) \,\cong\, \Z_2 \,$
and hence
\begin{equation}\label{equ:4.6}
  N^{\#}(\widetilde{f}_1, \widetilde{f}_2) \,\equiv\,
  \widetilde{N}(\widetilde{f}_1, \widetilde{f}_2) \,\equiv\,
  N(\widetilde{f}_1, \widetilde{f}_2) \,\not\equiv\, 0.
\end{equation}

This whole discussion and, in particular, formulas \eqref{equ:4.3}--\eqref{equ:4.6} establish the claims in example \ref{exa:1.27} for
$\, m \,\leq\, 5 \,$
and
$\, m \,=\, 9 $.
The interested reader is encouraged to carry out the necessary computations of
$\, E^{\infty} \,$
and
$\, \Gamma \,$
in the remaining cases
$\, m \,=\, 6, 7, 9 \,$
as an exercise (use the results in Toda's book [T]	
and especially also the exact sequences (2.11) and (4.4)\,).
\hfill$\Box$\\

Formulas \eqref{equ:4.3}--\eqref{equ:4.5} prove also the claim of corollary \ref{cor:1.24} for
$\, n = 2 $.
Thus it remains to consider the case where
$\, m, n \geq 3 $.

Given
$ [\widetilde{f}] \,\in\, \pi_m(S^{d(n'+1)-1}) $,
put
\begin{equation*}
  [f] \;\;{:=}\;\; [p \comp \widetilde{f}] \;\in\; \pi_m(\KP(n') ), \;\;\; m,\, n' \geq 1.
\end{equation*}
When
$ (\K,\, n') \neq (\H,\, 1) $
the claims of theorem \ref{thm:1.18}(ii) hold for the pair
$ (\,*\,, f) $
without any further restrictions concerning
$ m, n' \text{ or } \K = \R,\, \C \text{ or } \H $.
Put
\begin{equation*}
  q \,=\, d(n' + 1) - 1 \,=\, n + d - 1
\end{equation*}
for short.

In order to prove claim a.) in corollary \ref{cor:1.24}, consider first the Whitehead square
\begin{equation*}
  [\widetilde{f}] \,=\, [\iota_q, \iota_q] \,\in\, \pi_{2q-1}(S^q).
\end{equation*}
If
$\, q \,$
is odd, then
$\, [\widetilde{f}] \,$
lies in the kernel of
$\, E^{\infty} \,$
and of the
$ \Z$--valued
homomorphism
$\, \gamma_2 \,$
and hence of
$\, \Gamma \,$
(since
$\, 2[\iota_q, \iota_q] = 0 \,$,
cf. [Wh], p. 474 and 485);	
if, in addition,
$\, q \,\neq\,1, 3, 7 \,$
then
$\, [\iota_q, \iota_q] \,\neq\, 0 \,$
(by the famous result of F. Adams on odd Hopf invariants and an EHP--sequence argument); thus
\begin{equation*}
  N^{\#}(*, f) \,=\, R \cdot N^{\#}(\widetilde{*}, \widetilde{f}) \;\neq\; 0 \,=\, \widetilde{N}(*, f) \,=\, \widetilde{N}(\widetilde{*}, \widetilde{f}).
\end{equation*}
The remainder of claim a.) follows also from \ref{exa:1.15} and [K2], 1.17.	

Next let
$\, \widetilde{f} \,\colon\, S^{q+1} \longrightarrow S^q \,$
(and
$\, \widetilde{f} \,\colon\, S^{q+3} \longrightarrow S^q \,$,
resp.) be the iterated suspension of the complex Hopf map
$\, h_{\C} \,\colon\, S^3 \longrightarrow S^2 \,$
if
$\, \K = \R \,$
or if
$\, \K = \H \,$
and
$\, n' > 1 \,$
(and the composite of three such suspensions if
$\, \K = \C \,$,
resp.).
Then
$\, \widetilde{N}(*, f) \,=\, R \cdot \widetilde{N}(\widetilde{*}, \widetilde{f}) \,\neq\, 0 \,$
(since
$\, E^{\infty}([\widetilde{f}]) \,=\, \eta \,$
(or
$\, \eta^3 \,$,
resp.) and hence
$\, \Gamma([\widetilde{f}]) \,$
do not vanish), but
$\, N(*, f) \,=\, 0 \,$
(since
$\, 2\eta \,=\, 0 \,$
and
$\, \eta^4 \,=\, \nu\eta \,\in\, \pi_4^S \,=\, 0 \,$,
cf. \ref{thm:1.18} and [T]).	
If
$\, \K \,=\, \H,\; n' \,=\, 1 \,$
and hence
$\, \KP(n') \,=\, S^4 \,$,
choose
$\, [f] \,=\, 24 \cdot [h_{\H}] \,\in\, \pi_7(S^4) \,$;
then the Hopf invariant of
$\, [f] \,$
and hence
$\, \Gamma([f]) \,$
do not vanish but
$\, E^{\infty}([f]) \,\in\, \pi_3^S \,\cong\, \Z_{24} \,$
does; again
$\, \widetilde{N}(*, f) \,\neq\, N(*,f) \,$
(compare the proof of theorem \ref{thm:1.18}).	
This proves claim b.) in corollary \ref{cor:1.24}.

Finally, according to Toda [T], p. 177, lines 20 -- 25, and lemma 13.5,	
there exists an element
$\, \alpha_1(3) \,\in\, \pi_6(S^3) \,\cong\, \Z_{12} \,$
such that
$\, E^{\infty}(\alpha_1(3)) \,\in\, \pi_3^S \,\cong\, \Z_{24} \,$
has order
$ 3 $.
Choose
$\, [\widetilde{f}] \,$
(or
$\, [f] \text{ if } \K \,=\, \H,\, n = 4 \,$)
to be a suitable suspension of
$\, \alpha_1(3) \text{ if } \K = \R \,$
and of the Hopf classes
$\, [h_{\C}] \,\in\, \pi_3(S^2) \,$
and
$\, [h_{\H}] \,\in\, \pi_7(S^4) \,$
if
$\, \K \,=\, \C \text{ or } \H \,$,
resp.
Then
$\, N(*, f) \,>\, N^{\Z}(*,f) \,=\, 0 \,$
(since
$\, 2E^{\infty}(\alpha_1(3)),\, \eta^2,\, \nu^2,\, \nu \,\in\, \pi_*^S \,$
are nontrivial; cf. \ref{thm:1.18}(ii), IV and V).
This completes the proof of corollary \ref{cor:1.24}.
\hfill$\Box$\\

%
%
%

\section*{Acknowledgements}
This research was supported in part by DFG and DAAD.\\
Also it is a pleasure to thank M. Crabb, K. Knapp and D. Randall for helpful references.

The papers [K1]--[K6] and [KR] are available online at\\ http://www.uni-siegen.de/fb6/rmi/topologie/publications.html\,.

\newpage
\begin{thebibliography}{Koschorke}
\section*{References}

\bibitem[B1]{b1} \textit{R. Brown},
	Wecken properties for manifolds,
	Contemp. Math. \textbf{152} (1993), 9-21.

\bibitem[B2]{b2} ----------,
	Nielsen fixed point theory on manifolds, Banach Center Publ. \textbf{49} (1999), 19-27.

\bibitem[BS]{bs} \textit{R. Brown} and \textit{H. Schirmer},
	Nielsen coincidence theory and coincidence-producing maps for manifolds with boundary,
	Topology Appl. \textbf{46} (1992), 65-79.

\bibitem[Br]{br} \textit{R. Brooks},
	On removing coincidences of two maps when only one, rather than both, of them may be deformed by a homotopy,
	Pacific J. Math. \textbf{39}, no.3 (1971), 45-52.

\bibitem[BGZ]{bgz} \textit{S. A. Bogatyi}, \textit{D. L. Gon\c{c}alves} and \textit{H. Zieschang},
	Coincidence theory: the minimizing problem,
	Proc. Steklov Inst. Math. \textbf{225} (1999), 45-77.

\bibitem[C]{c} \textit{M. Crabb},
	The homotopy coincidence index,
	J. Fixed Point Theory Appl. \textbf{7} (2010), 1-32.

\bibitem[D]{d} \textit{J. P. Dax},
	Etude homotopique des espaces de plongements,
	Ann. Sc. Ec. Norm. Sup. \textbf{5} (1972), 303-377.

\bibitem[DG]{dg} \textit{A. Dold} and \textit{D. Gon\c{c}alves},
	Self-coincidence of fibre maps, Osaka J. Math. \textbf{42} (2005), 291-307.

\bibitem[GR 1]{gr1} \textit{D. Gon\c{c}alves} and \textit{D. Randall},
	Self-coincidence of maps from $S^q$-bundles over $S^n$ to $S^n$,
	Bol. Soc. Mat. Mexicana (3) \textbf{10} (2004), 181-192.

\bibitem[GR 2]{gr2} ----------,
	Self-coincidence of mappings between spheres and the strong Kervaire invariant one problem,
	C. R. Math. Acad. Sci. Paris \textbf{342} (2006), 511-513.

\bibitem[HQ]{hq} \textit{A. Hatcher} and \textit{F. Quinn},
	Bordism invariants of intersections of submanifolds,
	Trans. AMS \textbf{200} (1974), 327-344.

\bibitem[Ja]{ja} \textit{I. M. James},
	The topology of Stiefel manifolds,
	London Math. Soc. Lect. Notes \textbf{24}, Cambridge Univ. Press 1976.

\bibitem[Ji1]{ji1} \textit{B. Jiang},
	Fixed points and braids,
	Invent. Math. \textbf{75} (1984), 69-74.

\bibitem[Ji2]{ji2} ----------,
	Fixed points and braids II,
	Math. Ann. \textbf{272} (1985), 249-256.

\bibitem[Ji3]{ji3} ----------,
	Commutativity and Wecken properties for fixed points of surfaces and 3-manifolds,
	Topology Appl. \textbf{53} (1993), 221-228.

\bibitem[Ke]{ke} \textit{M. Kelly},
	Minimizing the number of fixed points for self-maps of compact surfaces,
	Pacific J. Math. \textbf{126} (1987), 81-123.

\bibitem[K1]{k1} \textit{U. Koschorke},
	Selfcoincidences in higher codimensions,
	J. reine und angew. Math. \textbf{576} (2004), 1-10.

\bibitem[K2]{k2} ----------,
	Nielsen coincidence theory in arbitrary codimensions,
	J. reine und angew. Math. \textbf{598} (2006), 211-236.

\bibitem[K3]{k3} ----------,
	Nonstabilized Nielsen coincidence invariants and Hopf-Ganea homomorphisms,
	Geometry and Topology \textbf{10} (2006), 619-666.

\bibitem[K4]{k4} ----------,
	Selfcoincidences and roots in Nielsen theory,
	J. Fixed Point Theory Appl. \textbf{2} (2007), 241-259.

\bibitem[K5]{k5} ----------,
	Minimizing coincidence numbers of maps into projective spaces,
	Geometry \& Topology Monographs \textbf{14} (2008), 373-391.

\bibitem[K6]{k6} ----------,
	Minimum numbers and Wecken theorems in topological coincidence theory.I,
	J. Fixed Point Theory Appl. \textbf{10} (2011), 3-36.

\bibitem[KR]{kr} \textit{U. Koschorke} and \textit{D. Randall},
	Kervaire invariants and selfcoincidences,
	Oberwolfach (2007) and New Orleans (2011), to appear.

\bibitem[Kn]{kn} \textit{K. H. Knapp},
	Vektorb\"undel--eine Einf\"uhrung, in preparation,

\bibitem[N]{n} \textit{J. Nielsen},
	Untersuchungen zur Topologie der geschlossenen zweiseitigen Fl\"achen,
	Acta Math. \textbf{50} (1927), 189-358.

\bibitem[S]{s} \textit{H. A. Salomonsen},
	Bordism and geometric dimension,
	Math. Scand. \textbf{32} (1973), 87-111.

\bibitem[T]{t} \textit{H. Toda},
	Composition methods in homotopy groups of spheres,
	Annals of Mathematics Studies 49, Princeton University Press 1962.

\bibitem[We]{we} \textit{F. Wecken},
	Fixpunktklassen. I, II, III,
	Math. Ann. \textbf{117} (1941), 659-671; \textbf{118} (1942), 216-234 and 544-577.

\bibitem[Wh]{wh} \textit{G. Whitehead},
	Elements of homotopy theory,
	Graduate Texts in Mathematics, Springer-Verlag, 1978.

\bibitem[Z]{z} \textit{X. Zhang},
	The least number of fixed points can be arbitrarily larger than the Nielsen number,
	Acta Sci. Nat. Univ. Pekin. \textbf{1986}, 15-25.


\end{thebibliography}
\end{document}